\documentclass[12pt,letter]{article}
\usepackage{amssymb,amsfonts,amsmath,setspace} 

\newcommand\proof{{\em Proof}}

\renewcommand\le{\leqslant}
\renewcommand\ge{\geqslant}

\newcounter{Remark}

\newtheorem{Theorem}{Theorem}
\newtheorem{Lemma}{Lemma}
\newtheorem{Remark}{Remark}

\numberwithin{equation}{section}
\numberwithin{Theorem}{section}
\numberwithin{Lemma}{section}
\numberwithin{Remark}{section}

\newcommand{\f}{\varphi}
\newcommand{\si}{\sigma}
\newcommand{\D}{\Delta}                  
\newcommand{\ind}{{\bf 1}}
\newcommand{\exn}{{\bf E}\,}
\newcommand{\pr}{{\bf P}}
\newcommand{\R}{{\mathbb{R}}}
\newcommand{\Z}{{\mathbb{Z}}}
\newcommand{\deq}{\stackrel{d}{=}}

\newcommand{\lez}{\preccurlyeq}
\newcommand{\grz}{\succcurlyeq}

\newcommand{\m}{\medskip}

\newcommand{\ati}{{\boldmath$[ \mbox{\boldmath$\psi $} ]$}}
\newcommand{\ado}{{\boldmath$[ \mbox{\boldmath$\psi,\downarrow$}   ]$}}
\newcommand{\aup}{{\boldmath$[ \mbox{\boldmath$\psi,\uparrow$}   ]$}}
\newcommand{\acu}{{\boldmath$[ \mbox{\boldmath$\psi,\uparrow \cup\downarrow$}   ]$}}

\newcommand{\Roh}{{\boldmath$[ {\mathbf R}$}$_{\alpha, \rho}${\boldmath$]$}}

\begin{document}
%
%
%
%

\title{Blackwell-type Theorems for\\ Weighted Renewal Functions$^1$}

\author{Alexander A.\ Borovkov$^2$ and  Konstantin A.\ Borovkov$^3$}

\date{}
\maketitle

\footnotetext[1]{Research supported by the President of the Russian Federation Grant NSh-3695.2008.1,  the  Russian Foundation for Fundamental Research Grant 08--01--00962 and the ARC
Centre of Excellence for Mathematics and Statistics of Complex Systems.}

\footnotetext[2]{Sobolev Institute of Mathematics,  Ac.~Koptyug avenue~4, 630090 Novosibirsk, Russia. E-mail: borovkov@math.nsc.ru.}

\footnotetext[3]{Department of Mathematics and Statistics, The University of Melbourne, Parkville 3010, Australia. E-mail: borovkov@unimelb.edu.au.}

\begin{abstract}
For a numerical sequence $\{a_n\}$ satisfying  broad assumptions on its ``behaviour on average'' and a random walk $S_n=\xi_1 +\cdots +\xi_n$ with i.i.d.\ jumps $\xi_j$ with positive mean $\mu$, we establish the asymptotic behaviour of the sums
\[
\sum_{n\ge 1} a_n \pr (S_n\in[x, x+\D))  \quad \mbox{as}\quad x\to \infty,
\]
where~$\D>0$ is fixed.
The novelty of our results is not only in much broader conditions on the weights $\{a_n\}$, but also in that neither the jumps~$\xi_j$ nor the weights~$a_j$ need to be positive. The key tools in the proofs are integro-local limit theorems and large deviation bounds.
For the jump distribution~$F$, we consider conditions of four types:  (a)~the second moment of~$\xi_j$ is finite, (b)~$F$ belongs to the domain of attraction of a stable law,   (c)~the tails of  $F$ belong to the class of the so-called locally regularly varying functions, (d)~$F$ satisfies the moment Cram\'er condition. Regarding the weights, in cases  (a)--(c) we assume that  $\{a_n\}$ is a so-called
$\psi$-locally constant on average sequence,   $\psi(n)$ being the scaling factor ensuring convergence of the distributions of $(S_n - \mu n)/\psi (n)$ to the respective stable law. In case~(d) we consider sequences of weights of the form  $a_n=b_n e^{qn},$  where $\{b_n\}$ has the properties assumed about the sequence $\{a_n\}$ in cases (a)--(c) for~$\psi(n)=\sqrt{n}.$


\end{abstract}

\section{Introduction}

Consider the following motivating problem where we will meet the objects to be studied in the present paper.  Let $\{\xi_j\}_{j\ge 1}$ and $\{\tau_j
\}_{j\ge 0}$ be independent   sequences  of random variables (r.v.'s), $\xi_j\deq \xi$ being independent identically distributed with a common distribution function $F$ and finite mean $\mu:= \exn \xi >0$, $\tau_j>0$ being arbitrary with  $a_j:=\exn \tau_j <\infty.$ Set
\[
S_n := \sum_{j=1}^n \xi_j, \qquad T_n := \sum_{j=0}^n \tau_j,\qquad n=0,1,2,\ldots,
\]
and consider the generalised renewal process  $S_{\nu (t)}$, where
$\nu (t):= \inf\{n\ge 0: T_n >t\}$, $t\ge 0$. Clearly, the process describes the movement of a particle that rests at zero during the time interval  $[0,T_0),$ then changes its location at time $T_0$ by a jump of size  $\xi_1$ and stays at the new location during the time interval  $[T_0, T_1),$ then jumps by $\xi_2,$ and so on. Processes of that kind are often encountered in applications, e.g.\ in queueing.

The time spent by the process  $S_{\nu (t)}$ in a given interval
\[
\D [x):= [x, x+\D), \qquad x\in \R, \quad \D>0,
\]
is equal to
\[
\tau (\D [x)) :=\sum_{n=0}^\infty \tau_n \ind (S_n \in \D[x)),
\]
where $\ind (A)$ is the indicator of the event~$A$. It is often important to know the behaviour of the expectation
\[
h(x, \D):= \exn  \tau (\D [x))= \sum_{n=0}^\infty a_n \pr (S_n \in \D[x))
\]
of that time and, in particular, its asymptotics as~$x\to\infty.$

Observe also that, if we ``swap time with space'' and assume that $\xi \ge 0$, then  $\tau (\D [x))$ will be the increments of the generalised renewal process  $T_{\eta (x)}$ on the intervals $\D [x)$, $\eta (x)$ being the renewal process generated by the sequence~$\{\xi_j\}$.  In that case, the assumption $\tau_j >0$ (so that~$a_j>0$) is not needed. In what follows, we will be allowing negative values for~$a_n.$

When the series
\[
H(x) := \sum_{n=0}^\infty a_n \pr (S_n <x)
\]
converges,  its sum $H(x)$ is referred to as the weighted renewal function for the sequence~$\{\xi_j\}$. In that case,
$h(x,\D)=H(x+\D) - H(x)$ is just the increment of the  function~$H$ on the interval~$\D[x)$. Note, however, that $h(x,\D)$ can exist even when $H(x)$ does not. Indeed,  if the weights~$a_n$ are bounded then a sufficient condition for $H(x)$ to exist and be finite is that  $\exn (\xi^-)^2 <\infty$, where $v^-:=\min\{v, 0\},$ while if $a_n\equiv 1$ and $\exn (\xi^-)^2=\infty$, then  $H(x)=\infty$ (see e.g.\ \S\,10.1 in~\cite{BoU}). On the other hand, the  quantity  $h(x,\D)$ is always finite provided that $0<
\mu <\infty,$ $|a_n|\le C <\infty$ (this follows from the transience of the random walk $\{S_n\}$, see e.g.\ Section~VI.10 in~\cite{Fe}).

There is substantial literature devoted to studying the asymptotics of~$H(x)$ and $h(x,\D)$ as $x\to\infty$,   recent surveys of the area and some further references being available in~\cite{OmTe,Li}. In what follows, we will mostly be  mentioning results referring to~$h(x, \Delta)$.

In the special case $a_n=1/n$, $H(x)$ is called the harmonic renewal function. It is closely related to the concepts of  factorization identities, ladder epochs and heights etc. That case was dealt with in~\cite{GrOmTe82, Gr88, Al91, BaOm01}.

In the case of regularly varying function (r.v.f.)  $a_x=x^{\gamma} L(x)$ (here $L$ is slowly varying), it was shown in~\cite{Al92} that, if $\gamma \ge 0$ and $a_x$ is ultimately increasing, then the condition
\begin{equation}
\label{conda}
\exn (\xi^-)^2 a_{\xi^-} <\infty
\end{equation}
is necessary and sufficient for having $H(x)<\infty$ for all $x$ (thus extending the above-mentioned results on the finiteness of $H$) which, in its turn, is equivalent to the asymptotics
 \begin{equation}
\label{asympa}
H(x)\sim \frac{x }{\mu (\gamma+1)}\, a_{x/\mu}, \qquad x\to\infty,
\end{equation}
while if $\gamma \in (-1, 0]$ and $a_x$ is ultimately decreasing, relation~\eqref{conda} still implies~\eqref{asympa}.
Here and in what follows, $\sim$ denotes asymptotic equivalence: we write  $f(x)\sim g(x)$  if $f(x)/g(x)\to 1$ as $x\to\infty$. We will also use the convention that $b_x:=b_{\lfloor x\rfloor}$ for a sequence $\{b_n\}$, where
$\lfloor x\rfloor$ is the integer part of~$x$.

For the function $h (x,\D)$, it was shown in~\cite{Al92} that, in the non-lattice case, if $a_x$ is an ultimately increasing r.v.f.\ of index  $\gamma \ge 0$, then condition
\begin{equation}
\label{condaa}
\exn  \xi^- a_{\xi^-}  <\infty
\end{equation}
is equivalent to the asymptotics
\begin{equation}
\label{lim}
h(x ,\D) \sim \frac{\D}{\mu}\, a_{x/\mu}, \qquad x\to\infty
\end{equation}
(with an analog of the relation holding true in the arithmetic case), and that \eqref{lim} also holds when   $\gamma \in (-1, 0]$ and $a_x$ is ultimately decreasing (in that case  \eqref{condaa} is  always met in the finite mean case). Earlier research for that case was done in~\cite{Ka72} (for~$a_x = x^{\gamma}$) and~\cite{MaOm84, EmMaOm84}.

Clearly, \eqref{lim} and its analog in the arithmetic case are direct extensions of the
celebrated Blackwell theorem in Renewal Theory that describes the asymptotics of $h(x
,\D)$ as  $x\to\infty$ in the case where~$a_n\equiv 1.$

The special case  $a_n = n!{n+k-1 \choose k}$ was investigated in~\cite{Sg97}, where
three-term approximations where obtained for $H(x)$ for nonlattice~$F$ with a finite
second moment.

In \cite{OmTe} it was shown that, in the non-lattice case, if the sequence $a_n$ is non-decreasing with $\lim_{n\to \infty}a_n = \infty$ and
\[
\lim_{s\downarrow 0} \limsup_{x\to\infty} \frac{a_{x(1+s)}}{a_{x}} = 1 ,
\]
then~\eqref{lim} holds true. That paper also showed that~\eqref{lim}  is true in  the
case where $a_n\downarrow 0$ as $n\to\infty$ (under some additional conditions).

When the sequence $a_n$ is fast changing, the asymptotics of $h (x, \D)$ will be different from~\eqref{lim}. The special case where $a_n =A^n$ for a fixed $A>0$,   $\pr(\xi \ge 0)=1$ and $F$ has an absolutely continuous component or is arithmetic, was considered in~\cite{Na68}. It was shown there, in particular, that, provided that $F$ has finite moment generating function and $A$ belongs to a suitable range of values, the asymptotics of $h(x, \Delta)$ will also be of exponential nature (the precise form of the asymptotics is given in our Section~\ref{cramer} below, see  relations~\eqref{lim_cr},  \eqref{lim_ar} with $b_n\equiv 1$).

Of the known general results, the following assertion stated as Theorem~6.1
in~\cite{Li} is the closest one to the main topic of the present paper:  {\em If the
sequence  $\{a_n\}$ is regularly oscillating in the sense that }
\begin{equation}
\label{reg_os}
\lim_{x,y\to\infty,\ x/y\to 1}\frac{a_x}{a_y} =1
\end{equation}
{\em and, as  $
x\to\infty,$ one has
\begin{equation}
\label{lin}
F_+ (x) := \pr (\xi \ge x) = o \left( \frac{a_x}{A_x}\right),
\quad\mbox{where} \quad A_n:=\sum_{j=0}^n a_j,
\end{equation}
then, in the non-lattice case, for any fixed~$\D >0$ one has} \eqref{lim}.

A relation analogous to~\eqref{lim} holds under the same assumptions in the arithmetic
case as well (Theorem~3.1 in~\cite{Li}).  In fact, it were these results that drew our
attention to the problem on the asymptotic behaviour of~$h(x,\D)$ as $x\to\infty$, as
we realized that the conditions imposed in~\cite{Li} on the weight sequence $\{a_n\}$
could be substantially relaxed provided that $F$ belongs to the domain of attraction of
a stable law.

In the present paper, we establish asymptotics of the form~\eqref{lim} under much more
general conditions on $\{a_n\}$ than have previously been dealt with and allowing $\xi$
to assume values of both signs. In Section~\ref{normal} we consider two cases that can
be treated using the same approach based on the integro-local Gnedenko--Stone--Shepp
theorems and some large deviation bounds from~\cite{BoBoBoo}: the case of finite
variance, and the case where $F$ belongs to the domain of attraction of a stable law
with index  $\alpha \in (1,2)$. In these cases, using the fact that, owing to the
integro-local theorems, the sequence of probability values $\pr (S_n\in[x, x+\D))$
varies in a nice regular way, we demonstrate \eqref{lim} under the assumption that the
weight sequence is ``$\psi$-locally constant on average'' (see the definition of that
property  in Section~\ref{normal})  and is  non-decreasing or non-increasing ``on
average''. Note that, under the assumptions made with regard to $\{a_n\}$, the above
conditions on the distribution of~$\xi$ are close to the ones necessary
for~\eqref{lim}. In Section~\ref{locrvf} we obtain necessary and sufficient conditions
for~\eqref{lim} to hold under the assumption that $\exn \xi^2 <\infty$ and the tails of
$F$ are of locally regular variation (see the definition thereof in the beginning of
Section~\ref{locrvf}). In~Section~\ref{cramer} we find the asymptotics of $h(x,\D)$
under the assumptions that  $F$ satisfies the moment  Cram\'er condtion and $a_n=b_n
e^{qn}$, where $\{b_n\}$ is a $\psi$-locally constant sequence with~$\psi
(n)=\sqrt{n}$, $q=\mbox{const}\neq 0$.

\section{The case where the second moment is finite or there is convergence to a non-normal stable law}  \label{normal}

To formulate appropriate conditions on the weight sequence, first recall the definition of an asymptotically  $\psi$-locally constant ($\psi$-l.c.)
function (see Definition~1.2.7 in~\cite{BoBoBoo}). Let $\psi
(t)>1,$ $t>0,$ be a fixed non-decreasing function.

{\em A function  $g(x)>0$ is said to be  $\psi$-l.c.\ if, for any fixed  $v\in\R$ such that $x+v\psi (x) \ge cx$ for some  $c>0$ and all large enough~$x$, one has}
\[
\lim_{x\to \infty} \frac{g(x+ v\psi (x))}{g(x)} =1.
\]
A sequence $\{u_n\}$ is called  $\psi$-l.c.\ if
$g(x):=u_x$ is a $\psi$-l.c.\ function. Everywhere in what follows  $\psi$ will be assumed to be an r.v.f.

Note that $\psi$-l.c.\ functions with $\psi\equiv 1$ are sometime referred to as ``long-tailed''. Furthermore, it follows from Theorem~1 in~\cite{BoBoPap} that, under broad assumptions on~$\psi$, a $\psi$-l.c.\ function will, in the terminology of~\cite{FoKoZa}, be ``$h$-insensitive'' (or ``$h$-flat'') with $h\equiv\psi$  (see Definition~2.18 in~\cite{FoKoZa}).

Introduce the following conditions on the weight sequence.

\m

{\bf Condition} \ati\  is satisfied for  sequence  $\{a_n\}$ if {\em  there exists an r.v.f.\   $d (t)$ such that $d(t)  =o (\psi (t)) $ as
$t\to\infty$ and the ``averaged sequence''
\begin{equation}
\label{tia}
\tilde a_n := \frac1{d (n)}\sum_{n\le k< n+d (n)} a_k >0
\end{equation}
is $\psi$-l.c.} One can always assume that the averaging interval length $d(n)$ is integer-valued.

\m

Sequences satisfying this condition  we will  call ``$\psi$-locally constant  on average''.
In the present section, the function  $\psi$ will be chosen to be the scaling from the respective limit theorem for the sequence of partial sums~$S_n$.

Any $\psi$-l.c.\ sequence clearly is a $\psi$-l.c.\ on average sequence: it satisfies condition \ati\ with $d(n)\equiv 1$. A simple example of a $\psi$-l.c.\ on average sequence that is not $\psi$-l.c.\ is provided by a periodic sequence   $a_n =a_{n-\lfloor n/d\rfloor d}$ with a period
$d\ge 2,$ such that there are distinct values among  $a_0, a_1, \ldots, a_{d-1}.$ In that case,  $\tilde a_n=\mbox{const},$ and
$d(n)$ can be chosen to be a (constant) multiple of~$d$. Note that condition \ati\
does not exclude the case where some of the weights $a_n$ can be negative.

Sequences satisfying one of the following two conditions could be called ``monotone on average''.

\m

{\bf Condition~\ado} \  is met for   $\{a_n\}$ if that
sequence {\em satisfies  \ati\ and, for some $r>1 $
and $c<\infty$, one has $|a_k| \le  c\tilde a_n$ for all $k> n/r.$
}

\m

{\bf Condition~\aup} \  is met for   $\{a_n\}$  if that
sequence  {\em satisfies  \ati\ and, for some $r>1 $
and $c<\infty$, one has $|a_k| \le c\tilde a_n$ for all $k<   nr.$
}

\m

By \acu\ we will denote the condition that at least one of   conditions \ado\ and \aup\ is satisfied for the sequence in question.

We will also need the following conditions on the tails
\[
F_+ (t)  = \pr (\xi \ge t)  \quad \mbox{and} \quad F_-(t):=\pr (\xi <-t)
\]
of the jump distribution. We will assume that either
\begin{equation}
\sigma^2 := \exn (\xi -\mu)^2 <\infty
 \label{si_fin}
\end{equation}
or the next condition is met:

\m

{\bf Condition~\Roh.}\ \ {\em  The two-sided tail
\[
F_*(t) := F_-(t) + F_+ (t)
\]
is an r.v.f.\ with index $-\alpha$,
$\alpha\in (0,2),$ and there exists the limit}
\begin{equation*}
\lim_{t\to\infty} \frac{F_+ (t)}{F_*(t)}=: \frac{1}{2}(\rho+1)\in [0,1].
\label{RHO}
\end{equation*}

\m

Given that condition  \Roh\ is satisfied, we put
\begin{equation}
\label{b}
b(t) := \inf\{x>0: \, F^*(x) <1/t \}, \quad t>1
\end{equation}
(recall that $b(t)=t^{1/\alpha} l(t)$, where $l(t)$ is slowly varying as $t\to\infty$, see e.g.\ Theorem~1.1.4(v)
in~\cite{BoBoBoo}).

\m

Now let
\begin{equation}
 \label{Sqrtb}
\psi (t):=\left\{\begin{array}{ll}
\sigma \sqrt{t}  & \quad \mbox{if~\eqref{si_fin} is met},\\
b (t)  & \quad \mbox{if~\Roh\ is met}.
\end{array}
 \right.
\end{equation}
It is well-known that, given that either \eqref{si_fin} or
\Roh\ is satisfied, the sequence of the distributions of the scaled partial sums $(S_n-\mu n)/\psi (n)$ converges weakly as $n\to\infty$ to the respective stable law that we will denote by $\Phi$ (see e.g.\ Theorem~1.5.1 in~\cite{BoBoBoo}).

Moreover, the integro-local  Stone--Shepp theorem holds true as well:
{\em if $F$ is non-lattice then, for the approximation error $\epsilon_n
(x, \D)$  in  the representation
\begin{equation}
\label{StoneThm0}
 \pr (S_n \in \D[x) )
   = \frac{\D}{\psi (n)}\,  \phi
    \biggl(\frac{ x-\mu n}{\psi  (n)}\biggr)
    + \frac{\epsilon_n (x, \D)}{\psi  (n)}, 
\end{equation}
$\phi$ being the density of  $\Phi$, one has
\begin{equation}
\label{StoneThm}
 \lim_{n\to\infty} \sup_{\D\in [\D_1, \D_2]} \sup_{x}
  \big| \epsilon_n (x, \D) \big| =0
\end{equation}
for any fixed   $0<\D_1 <\D_2<\infty$ } (see e.g.\ Theorems~8.7.1 and 8.8.2 in~\cite{BoU}, or Theorem~6.1.2 in~\cite{BoBoBoo}).

In the arithmetic maximum span~1 case (where $\pr (\xi \in \Z)=1$ and g.c.d.$\{k_1-k_2:   \pr (\xi =k_1)\pr (\xi =k_2) >0\}=1$) an analog of the above theorem holds for probabilities   $\pr (S_n=k),$ $k\in \Z$ (Gnedenko's theorem, see e.g.\ Theorems~8.7.3 and 8.8.4 in~\cite{BoU}, or Theorem~6.1.1 in~\cite{BoBoBoo}).

Now we will state the main results of this section.

\begin{Theorem}
 \label{thm1}
Let the jump distribution  $F$ be non-lattice and condition \ado\ be satisfied
for~$\{a_n\}$. Furthermore, assume that one of the following two conditions is met:
either

{\rm (i)}  $\exn \xi^2 <\infty$ and the right tail $F_+$ of
$F$ admits a regularly varying majorant:
\begin{equation}
\label{rvf}
F_+ (t) \le V(t) := t^{-\alpha} L_V(t), \qquad t>0,
\end{equation}
where $\alpha >2$ and $ L_V (t)$ is slowly varying as $t\to\infty$,
such that
\begin{equation}
 \label{VaA}
 V(x)=o \biggl(\frac{\tilde a_{x}}{\widetilde A_x}\biggr) \quad
  \mbox{as\quad $x\to\infty,$ }\quad \widetilde A_n:= \sum_{k\le n} \tilde a_k,
\end{equation}
or

{\rm (ii)} one has
\begin{equation}
 \label{VaA+}
 \mbox{\Roh\ holds for some $\alpha \in (1,2)$, \quad }
 F_+ (x)=o \biggl(\frac{\tilde a_{x}}{\widetilde A_x}\biggr) \ \mbox{as \ $x\to\infty.$}
\end{equation}
Then, for any fixed  $0<\D_1<\D_2<\infty,$ relation 
 \begin{equation}
\label{lim_til}
h(x ,\D) \sim \frac{\D}{\mu}\, \tilde a_{x/\mu}, \qquad x\to\infty,
\end{equation}
holds uniformly in~$\D\in [\D_1,\D_2].$
\end{Theorem}

\begin{Remark}
{\em We will write $f(x)   \preccurlyeq g(x) $ as $x\to\infty$ if
$f(x)\le c g(x)$ for some $c>0$ and all large enough~$x$.
If one has
\begin{equation}
 A_x \succcurlyeq \overline A_x : = \sum_{k\le x} |a_k| \quad \mbox{as\quad $x\to\infty,$}
  \label{9'}
\end{equation}
(this condition is close to \ado\ and is always met when $a_k \ge 0$),  then one can replace~\eqref{VaA} with
\begin{equation}
 V(x)=o \biggl(\frac{\tilde a_{x}}{  A_x}\biggr) \quad \mbox{as\quad $x\to\infty.$}
  \label{9''}
\end{equation}
Indeed, if~\eqref{9'} holds then, for $y=x+ d(x)$, we have
\[
A_y \grz  \overline A_y \ge  \widetilde{\overline A}_x \ge \widetilde{ A}_x,
\quad\mbox{where}\quad \widetilde{\overline A}_x:= \sum_{n\le x}  \widetilde{\overline a}_n,
  \quad \widetilde{\overline a}_n := \frac{1}{d(n)} \sum_{n\le k <n+ d(n)} |a_k|.
\]
Therefore, assuming  \ati\ and \eqref{9''} satisfied, one has
\[
V(x)   \preccurlyeq   V(y) =   o \biggl(\frac{\tilde a_{y}}{  A_y}\biggr)
  = o \biggl(\frac{\tilde a_{x}}{  A_y}\biggr)
   =  o \biggl(\frac{\tilde a_{x}}{  \widetilde A_x}\biggr),
\]
where the first relation holds since   $V$ is an r.v.f.\ and  $x\sim y$, and the third
one follows from the fact that   $\{\tilde a_x\}$ is a $\psi$-l.c. Thus~\eqref{VaA} is
established. It is not hard to see that condition $V(x) =  o \bigl( {\tilde a_{x}}/{
\overline A_x}\bigr)$ is also sufficient for~\eqref{VaA}.
 }
\end{Remark}

\begin{Theorem}
 \label{thm2}
Let the jump distribution $F$ be non-lattice and condition \aup\ be satisfied
for~$\{a_n\}$. Furthermore,  assume that one of the following two conditions is met:
either

{\em (i)} $\exn \xi^2 <\infty$ and the left tail  $F_-$ of $F$ admits a regularly varying majorant:
\begin{equation}
\label{rvf1}
F_-(t) \le W(t) := t^{-\beta} L_W(t), \qquad t>0,
\end{equation}
where $\beta >2$ and $L_W(t)$ is slowly varying as $t\to\infty$, such that
\begin{equation}
 \label{aW}
 \sum_{n\ge 0}\tilde a_n W(n) <\infty,
\end{equation}
or

{\rm (ii)}  one has
\begin{equation}
 \label{aW+}
 \mbox{\Roh\ holds for some $\alpha \in (1,2)$, \quad }  \sum_{n\ge 0}\tilde a_n F_-(n) <\infty.
\end{equation}
Then, for any fixed  $0<\D_1<\D_2<\infty,$ relation  \eqref{lim_til} holds uniformly in~$\D\in [\D_1,\D_2].$
\end{Theorem}

\begin{Remark}{\em If   $a_n\ge 0$ then, in  conditions~\eqref{aW} and \eqref{aW+}, one can replace  $\tilde a_n$ with~$a_n$. Indeed, since~$d(t)$ and~$W(t)$ are r.v.f.'s, adding up the coefficients of~$a_n$ in relation~\eqref{aW} we obtain
\[
\sum_{n\ge 0} \tilde a_n W(n) = \sum_{n\ge 0} a_n e_n W(n) ,
\]
where $e_n\to 1$ as $n\to\infty$. It follows from here that the series  $\sum  \tilde a_n W(n)$ converges iff $\sum a_n e_n W(n)$ does. Also, it is not hard to see that  convergence of $\sum |a_n| W(n)$ is sufficient for~\eqref{aW}. Similar remarks apply to condition~\eqref{aW+} as well.
}
\end{Remark}

\begin{Remark}{\em If, instead of assuming that $F$ is non-lattice, one assumes that $F$ is arithmetic, then the assertions of Theorems~\ref{thm1} and~\ref{thm2} remain true provided that both  $x$ and $\D$ are integer. In other words, in the arithmetic case, for integer-valued $x\to\infty$ one will have
\[
\sum_{n=0}^\infty a_n \pr (S_n=x) \sim \frac{1}{\mu}\, \tilde a_{x/\mu}.
\]
}
\end{Remark}

\proof\ {\em of Theorems~\ref{thm1} and~\ref{thm2}.} Put
\begin{equation}
\label{mn}
 n_\pm := \frac{x}{\mu}  \pm N \psi (x), \qquad m_\pm   := \frac{x}{\mu}\,r^{\pm 1},
\end{equation}
where $N=N(x)\to\infty$ slowly enough as $x\to\infty$  (the choice of
$N$ will be discussed below), $r>1$ is the quantity from conditions~\ado,
\aup.

Represent $h(x, \D)$  as
\begin{equation}
 \label{u-16'}
h(x, \D)   =  \Sigma_{2-} + \Sigma_{1-} +\Sigma_{0} + \Sigma_{1+} + \Sigma_{2+},
\end{equation}
where the terms on the right-hand side of  \eqref{u-16'} are the sub-sums of the products $a_n \pr (S_n \in \D[x))$ over the following ranges of $n$ values:
\[
 \Sigma_{2-}:=\sum_{n<m_-}  , \
 \Sigma_{1-}:= \sum_{m_-\le n < n_-} , \
 \Sigma_{0}:= \sum_{n_-\le n < n_+}, \
 \Sigma_{1+}:= \sum_{n_+\le n < m_+}, \
 \Sigma_{2+}:= \sum_{n \ge m_+}.
\]
Lemmata~\ref{lem1}--\ref{lem5} are devoted to evaluating these subsums.  The assertions of Theorems~\ref{thm1} and~\ref{thm2} will follow from these results.

\begin{Lemma}
 \label{lem1}
If condition~\ati\ is satisfied for the sequence  $\{a_n\}$ and either condition~\eqref{si_fin} or condition~\Roh\ with $\alpha\in (1,2)$ is met then
\[
\Sigma_0 \sim \frac{\D}{\mu}\, \tilde a_{x/\mu},
\]
provided that the value  $N=N(x)$ in~\eqref{mn} tends to infinity slowly enough as~$x\to\infty.$
\end{Lemma}

\proof. Assume for simplicity that   $n_-$ is integer (as we will see from what follows, changing the values of $n_\pm$ and $m_\pm$ by amounts of the order $o(\psi (x))$ does not change anything in the proof).

Construct a  sequence  $n_k,$ $ 0\le k\le  K+1,$ by setting $n_0 := n_-$,
\[
n_{k+1}:= n_k + d(n_k), \quad k=0,1,\ldots, K:=\min\{k\ge 0: \, n_{k+1}\ge n_+\}
\]
(so that  $K\sim 2N\psi (x)/d(x/\mu)$), and amend somewhat  $n_+$
by putting, in agreement with the above remark, its value equal to $n_+:= n_{K+1}.$ Partition the set $[n_-, n_+)$ into semi-intervals
$[ n_k, n_{k+1}),$ $k=0, 1, \ldots, K.$ On each of these intervals, the probabilities
\[
p_n:= \pr (S_n \in \D[x))
\]
remain ``almost constant''   (in the ratio sense).  More precisely, putting
\[
\pi(n) :=  \phi \bigl( (x-\mu n)/\psi (n) \bigr)/\psi (n),
\]
we obtain by virtue of \eqref{StoneThm0}, \eqref{StoneThm} and the relation
$d(n_k) = o(\psi (x/\mu))$ that, for $n\in [n_k,n_{k+1}),$ one has
 \[
p_n =   (1+o(1)) p_{n_k} = (1+o(1))  \D \pi (n_k)
\]
uniformly in $k\in [0,K] $ provided that $N\to \infty$ slowly enough
as $x\to\infty.$ Hence the sub-sums $\sum_{n\in [n_{k}, n_{k+1})} a_n
p_n$,\ $k=0,1,\ldots, K,$\  that comprise  $\Sigma_0$ are of the form
\[
(1+o(1)) (n_{k+1}- n_k)\tilde a_{n_k}  \D \pi (n_k).
\]
But it follows from condition~\ati\ that  $\tilde a_{n_k} = (1+o(1)) \tilde
a_{x/\mu}$ uniformly in $k\le K$ provided that $N\to \infty$ slowly enough as $x\to\infty$ (see Theorem~1 in~\cite{BoBoPap}), so that
\[
\sum_{n\in [n_{k}, n_{k+1})} a_n p_n = (1+o(1)) \D \tilde a_{x/\mu} (n_{k+1}- n_k) \pi (n_k),
\]
and therefore
\begin{align}
\label{14'}
 \Sigma_0 & = (1+o(1)) \D \tilde a_{x/\mu}\sum_{k=0}^K  (n_{k+1}- n_k) \pi (n_k)
 \notag \\
 & = (1+o(1)) \D \tilde a_{x/\mu}\sum_{k=0}^K
   \frac{n_{k+1}- n_k}{\psi  (n_k)}\, \phi \biggl(\frac{ x-\mu n_k}{\psi  (n_k)}\biggr).
\end{align}
As
\[
\frac{n_{k+1}- n_k}{\psi  (n_k)}\sim \frac{n_{k+1}- n_k}{\psi  (x/\mu)}\to 0,
 \qquad
 \phi \biggl(\frac{ x-\mu n_k}{\psi  (n_k)}\biggr) \sim \phi \biggl(\frac{ x-\mu n_k}{\psi  (x/\mu)}\biggr)
 \]
uniformly in $k\le K$ (if $N\to \infty$ slowly enough  as
$x\to\infty$), the last sum in~\eqref{14'} is (up to the factor $(1+ o(1))$) a Riemann sum  for the integral
\[
\int_{x/\mu - N\psi (x)}^{x/\mu + N\psi (x)}
  \phi \biggl(\frac{ x-\mu t}{\psi  (x/\mu)}\biggr) \frac{dt}{\psi  (x/\mu)}
 =\frac{1+o(1)}{\mu}   \int_{- Nc}^{Nc}
  \phi (u) \, du
  \to \frac{1}{\mu} ,
\]
where $c=\lim_{x\to\infty} \psi (x) /\psi (x/\mu)= \mu^{1/\alpha}$
($\alpha=2$ in the case where $\exn \xi^2 <\infty$). The lemma  is proved.
\hfill$\square$

\begin{Lemma}
 \label{lem2}
Under the conditions of Lemma~\ref{lem1}, if\/~\acu\ is met then
\[
\Sigma_{1\pm}= o(  \tilde a_{x/\mu}), \qquad x\to\infty.
\]
\end{Lemma}

\proof. By virtue of \acu\ one has
\[
|\Sigma_{1\pm}|\lez   \tilde a_{x/\mu}  h_{\pm} (x),
 \quad \mbox{where} \quad h_- (x) := \sum_{m_-\le n < n_-} p_n,\quad
h_+ (x) := \sum_{n_+\le n < m_+} p_n.
\]
Clearly, $h_\pm (x)\le h(x) -h_0(x),$ where, as $x\to\infty$,
\begin{equation}
 \label{g}
h(x):= \sum_{n\ge 0} p_n \to \frac{\D}{\mu}
\end{equation}
by Blackwell's theorem, and
\begin{equation}
 \label{g0}
h_0 (x):= \sum_{n_- \le n< n_+} p_n \to \frac{\D}{\mu}
\end{equation}
by Lemma~\ref{lem1} for the sequence  $a_n\equiv 1.$
It follows from here that  $h_\pm (x) \to 0$ as  $x\to\infty.$ The lemma  is proved. \hfill$\square$

\m

To bound  $\Sigma_{2\pm}$ we will need the following extension of
Lemma~6.1 from~\cite{Li} to the case where $\pr (\xi <0) >0$. Set
\[
\underline{S}^{(n)} := \inf_{k\ge n} (S_k - S_n) \deq  \underline{S}^{(0)},
 \qquad
\overline{S}_n := \max_{k\le n} S_k.
\]
Since $\mu >0,$ the r.v.\ $\underline{S}^{(0)}$ is proper and
\begin{equation}
 \label{gamm}
\gamma:= \pr (\underline{S}^{(0)} =0)>0
\end{equation}
(see e.g.\ \S\,12.2 in~\cite{BoU}).

\begin{Lemma}
 \label{lem3}
If $\D>0$ is such that  $F_+ (\D) >0$ then, for any  $n\ge 1,$ one has
 \begin{align}
\label{ineq1}
\sum_{k\le n} \pr (S_k \in \D [x)) &\le \frac{\pr(\overline{S}_n \ge x)}{\gamma F_+(\D)},
 \\
 \sum_{k\ge n} \pr (S_k \in \D [x)) &\le
 \frac{\pr(S_n +\underline{S}^{(n)} < x+\D)}{\gamma F_+(\D)}.
\label{ineq2}
\end{align}
\end{Lemma}

\proof. For any $I\subset \R_+$ holds
\begin{align}
\sum_{k\in I} \pr (S_k \in \D [x))
    =  \Sigma' + \Sigma'',
   \label{SiSi}
\end{align}
where
\begin{align*}
\Sigma' & :=  \sum_{k\in I} \pr \bigl(S_k \in \D [x), \, S_{k+j}\not\in \D[x)\
   \mbox{\em for all}  \ j\ge 1\bigr),
  \\
\Sigma'' & := \sum_{k\in I} \pr \bigl(S_k \in \D [x),\, S_{k+j} \in
\D[x)\  \mbox{\em for some}  \ j\ge 1\bigr).
   \label{SiSi}
\end{align*}
Since the events in the probabilities in the sum $\Sigma'$ are mutually exclusive, one has
\[
\Sigma' \le \pr \biggl( \bigcup_{k\in I} \{S_k\in \D[x)\}\biggr)
 \le
 \left\{
 \begin{array}{ll}
 \pr (\overline{S}_n \ge x)  & \mbox{ if $I=(0,n]$},
 \vphantom{\displaystyle \sum_1}\\
 \pr(S_n +\underline{S}^{(n)} < x+\D)  & \mbox{ if $I=[n,\infty)$.}
 \end{array}
 \right.
\]
Further, as    $S_k$ and  $\xi_{k+1}+
\underline{S}^{(k+1)}  \deq \xi_{1} + \underline{S}^{(1)} $ are independent of each other, we obtain
\begin{align*}
\Sigma'' & \le  \sum_{k\in I}
 \pr \bigl(S_k \in \D [x),\, \inf_{j>0}S_{k+j}< x+\D\bigr)
 \\
 & \le
\sum_{k\in I} \pr (S_k \in \D [x),\, \xi_{k+1} + \underline{S}^{(k+1)} <\D)
  \\
  & = \pr (\xi_{1} + \underline{S}^{(1)} <\D) \sum_{k\in I} \pr (S_k \in \D [x)),
\end{align*}
where,  using definition~\eqref{gamm} and the independence of $\xi_{1} $ and $ \underline{S}^{(1)}  $, we have
\begin{align*}
\pr (\xi_{1} + \underline{S}^{(1)} <\D)
 & = \pr (\xi_{1} + \underline{S}^{(1)} <\D, \, \underline{S}^{(1)}=0)
  + \pr (\xi_{1} + \underline{S}^{(1)} <\D, \, \underline{S}^{(1)}< 0)
 \\
 & \le
 \pr (\xi_{1}  <\D, \, \underline{S}^{(1)}=0)
  + \pr (  \underline{S}^{(1)}< 0)\\
  & =
  \pr (\xi_{1}  <\D )\gamma + 1 - \gamma =  1 - \gamma F_+ (\D).
\end{align*}
Substituting  the bounds we established  for $\Sigma'$ and $\Sigma''$ into~\eqref{SiSi} first in the case where  $I=(0,n]$ and then where $I=[n,\infty)$, we obtain the assertion of the lemma.\hfill$\square$

\begin{Lemma}
 \label{lem4}
For the relation
\begin{equation}
 \label{Si2-}
\Sigma_{2-} = o( \tilde a_{x/\mu}), \qquad x\to\infty,
\end{equation}
to hold it suffices that one of the following conditions is met:

{\rm (i)}~$\{a_n\}$ satisfies \ado\ and one of conditions\/~{\rm(i), (ii)}  of Theorem~$\ref{thm1}$ is met;

{\rm (ii)}~$\{a_n\}$ satisfies \aup.
\end{Lemma}

\proof. (i)~Assume for simplicity (and without losing generality) that $x/\mu\equiv m_- r=r^{M+1}$   (see~\eqref{mn})  for some integer~$M\ge 1$, where $r>1$ is the quantity from condition~\ado, and consider the semi-intervals
\begin{equation}
 \label{Ij}
I_j := [r^{j-1}, r^j), \qquad j=1,2,\ldots
 \end{equation}

First let $\D>0$ be such that $F_+ (\D) >0$. Then, using condition \ado, inequality~\eqref{ineq1}  and notation $\chi_j:= \ind
(I_j\cap \mathbb{N} \neq \varnothing)$, we obtain that
\begin{align}
|\Sigma_{2-}| & \le \sum_{j=1}^M \sum_{n\in I_j} |a_n| p_n
 \lez  \sum_{j=1}^M \tilde a_{r^j} \sum_{n\in I_j}  p_n
 \notag \\
  & \le \sum_{j=1}^M \tilde a_{r^j} \pr (\overline{S}_{r^j} \ge x)
  \chi_j.
  \label{Si20}
\end{align}

Recall that, in the case of finite variance, if condition~\eqref{rvf} is satisfied and   $n\lez  y \to \infty$, then
\begin{equation}
 \label{414}
\pr \Bigl(\max_{k\le n} (S_k - \mu k)\ge y\Bigr)\lez  n V(y)
 \end{equation}
(see  Corollary~4.1.4(i) in~\cite{BoBoBoo}). Observing that
\[
\{\overline{S}_{r^j} \ge x\}
 = \{\overline{S}_{r^j} - \mu r^j \ge x - \mu r^j \}
  \subset  \Bigl\{\max_{k\le r^j} (S_k - \mu k) \ge x - \mu r^j\Bigr\}
 \]
and  $ x - \mu r^j\ge  x - \mu r^M=  (1- r^{-1})x$ for $j\le M$, we obtain from~\eqref{414} the bound
\[
\pr (\overline{S}_{r^j} \ge x)
 \le \pr \Bigl(\max_{n\le r^j}(S_n - \mu n) \ge (1- r^{-1})x\Bigr)
 \lez  r^j V(x),\quad j\le M.
\]
Substituting this bound into~\eqref{Si20} and using the inequality
$\tilde a_{r^j} \le c\tilde a_{n} ,$  $n\in I_j$ (it holds by virtue of condition~\ado), we find that
\[
|\Sigma_{2-}| \lez  V(x) \sum_{j=1}^M  r^j \tilde a_{r^j}
  \chi_j
 \lez  V(x/\mu) \sum_{j=1}^M  \sum_{n\in I_j}  \tilde a_{n}
  \le V(x/\mu) \tilde A_{x/\mu}.
\]
Now relation \eqref{Si2-} follows immediately from condition~\eqref{VaA}.

In the case of convergence to a non-normal stable law, the above argument remains valid provided that, when justifying  inequality~\eqref{414},  we replace the reference to Corollary~4.1.4 in~\cite{BoBoBoo}  with that to Corollary~3.1.2 from the same monograph.

If $F_+ (\D) =0$ (i.e.\ $\pr (\xi < \D) =1$) then one can always find a $k\in \mathbb{N}$ such that for  $\D_k:= \D/k$ one has $F_+ (\D_k) >0$. After that, it remains to apply the bound we have just derived to each of the terms on the right-hand side of the representation
\[
h (x, \D) = h (x, k\D_k) = \sum_{j=0}^{k-1} h (x + j\D_k, \D_k).
\]

(ii)~In this case,
\[
|\Sigma_{2-}| \le \sum_{n<m_-} |a_n| p_n\lez  \tilde a_{x/\mu}   \sum_{n<m_-}   p_n
 \le  \tilde a_{x/\mu} (h(x) - h_0 (x)) = o(\tilde a_{x/\mu})
\]
by virtue of \eqref{g} and \eqref{g0}. The lemma is proved.\hfill$\square$

\begin{Lemma}
 \label{lem5}
For the relation
\[
\Sigma_{2+} = o(\tilde a_{x/\mu}), \qquad x\to\infty,
\]
to hold it suffices that one of the following conditions is met:

{\rm (i)}~$\{a_n\}$ satisfies \ado;

{\rm (ii)}~$\{a_n\}$ satisfies \aup\ and one of conditions\/~{\rm(i), (ii)}  of Theorem~$\ref{thm2}$ is met.
\end{Lemma}

\proof. (i)~In this case,
\[
|\Sigma_{2+}| \lez  \tilde a_{x/\mu}\sum_{n>m_+} p_n
 \le \tilde a_{x/\mu} (h(x) - h_0(x)) =o(\tilde a_{x/\mu} )
\]
by virtue of~\eqref{g}, \eqref{g0}.

\m

(ii)~Using notation \eqref{Ij} and assuming for simplicity that now  $x/\mu\equiv m_+/ r=r^{M-2}$ for some integer~$M > 1$,
where $r>1$ is the quantity from condition~\aup, we obtain from condition~\aup\ and inequality~\eqref{ineq2} that
\begin{align}
|\Sigma_{2+}| & \le \sum_{j\ge M } \sum_{n\in I_j} |a_n| p_n
 \lez  \sum_{j \ge  M} \tilde a_{r^{j-1}}\sum_{n\in I_j}  p_n
 \notag
 \\
 & \lez  \sum_{j\ge M } \tilde a_{r^{j-1}}
 \pr \bigl(S_{r^{j-1}} +\underline{S}^{(r^{j-1})} < x+\D\bigr)
 \label{Si2+}
\end{align}
(in contrast to the proof of Lemma~\ref{lem4}, we do not need to use the indicators $\chi_j$ here since $M$ is large and therefore   $r^j - r^{j-1}>1$ for~$j\ge M$).

Note that, for $j\ge M$ and large enough $M$, the mean value $\exn S_{r^{j-1}} =\mu r^{j-1} \equiv x r$ of the first r.v.\ in the probabilities on the right-hand side of~\eqref{Si2+} exceeds the value $x+\D$ by
\[
  \mu r^{j-1} - (x+\D )
 =   \mu r^{j-1}   - \mu r^{M-2}  - \D
 \ge   \mu r^{j-2} (r-1) - \D
  \ge   2cr^j
\]
for some $c>0,$ whereas   $\underline{S}^{(r^{j-1})} $
is stochastically minorated by the global minimum~$\underline{S}^{0}$. Therefore, to bound the probabilities of the right-hand side of~\eqref{Si2+}, we will need to use large deviation results for the left tails of the distributions of the r.v.'s~$S_n$ and~$\underline{S}^{0}$. In the case of finite variance, the required bound for the former r.v.\ is contained in the assertion of Corollary~4.1.4(i) in~\cite{BoBoBoo} (applied to the random walk with jumps $-\xi_j$). To bound the left tail of $\underline{S}^{0}$, it suffices to use condition~\eqref{rvf1} to construct a random walk with i.i.d.\ jumps that have a positive mean and distribution with a regularly varying left tail, such that it stochastically minorates~$\{S_n\}$, and then make use of Theorem~7.5.1 in~\cite{BoBoBoo}.

Following this approach, we obtain the bound
\begin{align*}
 \pr \bigl(S_{r^{j-1}} +\underline{S}^{(r^{j-1})} < x+\D\bigr)
   & \le  \pr \bigl(S_{r^{j-1}} -\mu r^{j-1} <  - cr^j \bigr)
 +\pr \bigl( \underline{S}^{(0)} < -cr^j\bigr)\\
 & \lez
 r^{j-1} W(cr^{j-1}) + cr^{j} W(c r^{j}) \lez  r^{j} W(r^{j}).
\end{align*}

Now returning to~\eqref{Si2+}, we obtain
\[
|\Sigma_{2+}| \lez  \sum_{j\ge M } \tilde a_{r^{j-1}} r^{j} W(r^{j})
 \lez  \sum_{j\ge M } W(r^{j}) \sum_{n\in I_j}  \tilde a_n
 \lez  \sum_{n\ge m_+} \tilde a_n W(n),
\]
where we again used condition~\aup\ and the fact that  $W$ is an
r.v.f. Therefore, under condition~\eqref{aW}, one has
\[
\Sigma_{2+} =o(1) = o(\tilde a_{x/\mu}),
\]
where the last relation follows from~\aup.

In the case of convergence to a non-normal stable law, instead of Corollary~4.1.4 one should make use of Corollary~3.1.2 in~\cite{BoBoBoo}. The lemma is proved.  \hfill$\square$

\m

The assertions of Theorems~\ref{thm1} и \ref{thm2} follow from Lemmata~\ref{lem1}, \ref{lem2}, \ref{lem4} и~\ref{lem5}. The required uniformity in  $\D$ in the bounds established in Lemmata~\ref{lem2}--\ref{lem5} is obvious.

\section{The case where the   tails of the jump distribution  are of local regular variation}
\label{locrvf}

We will say that $V(x)$ is a  {\em locally regularly varying
function}  (l.r.v.f.) if it is an r.v.f.\ and, moreover,
for any fixed $\D>0$ holds
\begin{equation}
 \label{lrvf}
V(x ) - V(x +\D) = \D v(x) (1 + o(1)),\qquad  x\to\infty,
\end{equation}
where $v(x) = \alpha V(x)/x$, $-\alpha$ is the index of the  r.v.f.\ $V(x)$
(cf.~\eqref{rvf}). It will be assumed throughout this section that $\alpha >2$. Property~\eqref{lrvf} could be called the
``differentiability of  $V$ at infinity''.

It is clear that \eqref{lrvf} will hold if the slowly varying function $L_V$ in the representation on the right-hand side of~\eqref{rvf}
is differentiable and $L_V'(x) = o(L_V(x)/x)$ as $x\to\infty$.

If the tail $F_+$ (or $F_-$) of the distribution  $F$ is an l.r.v.f., then the derivation of the asymptotics of $h(x, \D)$  as $x\to\infty$ substantially simplifies due to the fact that integro-local theorems valid on the whole half-line are available for such distributions. For instance, if $F_+$
is an l.r.v.f.\ then, as it was established in~\cite{Mo}, in the case where $\exn \xi =0$ and $\sigma^2 =\exn\xi^2 <\infty$,  one has the following representation valid for $x> c\sqrt{n},$ $c=\mbox{const}>0:$
\[
\pr (S_n\in \D [x)) = \biggl[\frac{\D}{\si \sqrt{n}}\,
 \phi \biggl(\frac{x}{\si \sqrt{n}}\biggr) + n\D v(x)\biggr](1 + o(1)), \qquad n\to\infty.
\]
where $\phi$ is the standard normal density and
\[
v(x) = \frac{\alpha F_+ (x)}{x}
\]
is the  function from representation~\eqref{lrvf} for $V(x) = F_+ (x)$.
A similar assertion holds for  $\pr (S_n\in \D [x))$ in the range $x
\le - c\sqrt{n},$ provided that  $F_-$ is an l.r.v.f.; in that case, we will use notation
 \[
w(x) := \frac{ \beta  F_- (x)}{x},
\]
where $-\beta$ is the index of the r.v.f.~$F_-$.

From here and the Stone--Shepp theorem we immediately obtain the following assertions that reduces the problem of computing the asymptotics of
$h(x,\D)$ to a purely analytic exercise.

\begin{Theorem}
Let $\exn \xi =\mu >0,$ $\sigma^2 <\infty,$ and the tails $F_+ (x)=
x^{-\alpha} L_V (x),$ $F_- (x)= x^{-\beta} L_W (x)$ be
l.r.v.f.'s,   $\min\{\alpha, \beta \}>2$, $b(n):=\sigma \sqrt{n}$.
Then, for any fixed  $\D>0,$ $c_-\in (0,1/\mu) $ and
$c_+\in (1/\mu,\infty),$ one has
\begin{align}
h(x,\D)  & = \D\biggl[
 \frac{1}{\sigma \sqrt{2\pi n}} \sum_{c_- x\le n \le c_+ x} a_n e^{-(x-\mu n)^2/(2 \sigma^2 n)} \notag\\
  & +  \sum_{n> x/\mu + b(x)} a_n n w(\mu n-x)
  +  \sum_{n< x/\mu - b(x)} a_n n v(x- \mu n )
\biggr] (1+ o(1))
 \label{summa}
\end{align}
as $x\to\infty$.
\end{Theorem}

The next assertion shows that, in  Theorem~\ref{thm2}(ii), condition~\eqref{aW+}   on the relationship between the tails $F_-$ and weights  $a_n$ cannot be weakened, and that condition~\eqref{VaA+} on the relationship between  $a_n$ and $F_+$ can be extended to the minimal one under the assumptions of the present section. Since we assume that $F$ has a finite second moment, we automatically put $\psi (t) :=\sigma \sqrt{t}$ (cf.~\eqref{Sqrtb}).

\begin{Theorem}
\label{thmN}
Let $\exn \xi =\mu >0,$ $\sigma^2 <\infty,$ and $a_n\ge 0.$

{\rm (i)} If\/ $F_-$ is an l.r.v.f.\ then the divergence of the series   $\sum_n a_n F_-(n) =\infty$ implies that  $h(x,\D)=\infty$ for any~$\Delta >0.$

{\rm (ii)} If condition \aup\ is satisfied for $\{a_n\}$ and $F_-$ is an l.r.v.f.\ then \eqref{lim_til} holds iff the series $\sum_n a_n
F_-(n) $  converges.

{\rm (iii)} Let condition \ado\ be satisfied for $\{a_n\}$ and $F_+$ be an l.r.v.f. Then, as $x\to\infty$,
\begin{equation}
 \label{h+}
h(x,\D)=  \frac{\D}{\mu}\, \tilde{a}_{x/\mu} (1 + o(1))
 + \D \sum_{n< x/(r\mu)  } a_n n v(x- \mu n ) (1 + o(1)).
\end{equation}
In this case, relation  \eqref{lim_til} holds iff
\begin{equation}
 \label{F+o}
F_+ (x) = o\biggl(\frac{x\tilde{a}_{x/\mu}}{B_{x/(r\mu)}}\biggr),\qquad  x\to\infty,
\end{equation}
where $B_x :=\sum_{k\le x} k a_k.$
\end{Theorem}

\begin{Remark}
\label{RemB}   {\em Note that the principal part of the second term on the right-hand side of~\eqref{h+}  is always between the values
\begin{equation*}
\D  B_{x/(r\mu)} v(x) \quad\mbox{and}\quad \D  B_{x/(r\mu)} v((r-1)x).
\end{equation*}
}
\end{Remark}

\begin{Remark}
\label{RemC}  {\em  If the series  $A_\infty$ diverges and $a_x$ is an r.v.f.,  then conditions~\eqref{F+o} and \eqref{VaA} are equivalent. However, if $A_\infty < \infty$ then condition~\eqref{F+o} is broader than~\eqref{lin}, \eqref{VaA}. For example, if  $B_\infty<\infty$ then  $A_\infty<\infty,$
condition~\eqref{F+o} takes the form $F_+ (x) = o( x\tilde{a}_{x/\mu})$, whereas condition~\eqref{VaA} will mean that $F_+ (x) = o(\tilde{a}_{x/\mu})$.
 }
\end{Remark}


\proof{\ \em of Theorem~\ref{thmN}}. (i)~This assertion  follows from the  inequality
\[
h(x,\D) \ge \D \sum_{n\ge x/(r\mu)} a_n n w (\mu n - x) \grz
 \sum_{n\ge x/(r\mu)} a_n  F_- (  n  ).
\]

\medskip

(ii)~We will make use of the representation
\[
h(x,\D) = \Sigma_{1-} +  \Sigma_{0} + \Sigma_{1+} + \Sigma_{2+},
\]
where
\begin{align*}
\Sigma_{1-} &:= \sum_{ n   < x/\mu - Nb(x)} a_n p_n,&
\Sigma_{0}& := \sum_{|n - x/\mu| \le  Nb(x)} a_n p_n,\\
\Sigma_{1+}    & : =
    \sum_{ n\in ( x/\mu + Nb(x),\, x/\mu + c_1 \sqrt{x\ln x}]} a_n p_n, &
 \Sigma_{2+} & :=  \sum_{ n  > x/\mu + c_1 \sqrt{x\ln x}} a_n p_n.
\end{align*}
By Lemma~\ref{lem1}, provided that $N=N(x)\to\infty$ slowly enough as $x\to\infty$, one has
\[
 \Sigma_{0}
  \sim \frac{\D}{\mu}\, \tilde{a}_{x/\mu},
\]
whereas by Lemma~\ref{lem2}
\[
 \Sigma_{1-} 
 =o(\tilde{a}_{x/\mu}).
\]
Further, for  $c_1=\mbox{const} >0$ one has
\begin{align*}
 \Sigma_{1+}    &= \frac{\D (1+o(1))}{ \sigma \sqrt{2\pi n}}
    \sum_{ n\in ( x/\mu + Nb(x),\, x/\mu + c_1 \sqrt{x\ln x}]}  a_n e^{- (x-\mu n)^2/(2 \sigma^2 n)} \\
  & \lez   \tilde{a}_{x/\mu} e^{-c N^2}
  =o(\tilde{a}_{x/\mu})
\end{align*}
and
\begin{equation}
 \label{Sig2}
 \Sigma_{2+} = (1+o(1)) \D \sum_{ n  > x/\mu + c_1 \sqrt{x\ln x}}  a_n  n w(\mu n-x).
\end{equation}
If the series $\sum_{ n }  a_n  F_- (n)$ diverges then $\Sigma_{2+}
=\infty$ (see the proof of part~(i) above). Otherwise make use of the representation
\[
 \Sigma_{2+} =  \Sigma_{2,1+} +  \Sigma_{2,2+},
 \qquad  \Sigma_{2,1+}:= \sum_{n\in ( x/\mu + c_1\sqrt{x\ln x},\,  xr/\mu]},\qquad
   \Sigma_{2,2+}:= \sum_{n>  xr/\mu}.
\]
It is not hard to see that
\[
 \Sigma_{2,1+} \lez  \tilde{a}_{x/\mu}
 \sum_{n\in ( x/\mu + c_1\sqrt{x\ln x}, \, xr/\mu]}
 n w (\mu n -x) \lez   \tilde{a}_{x/\mu} x F_- (Nb(x))=o(\tilde{a}_{x/\mu})
\]
and, by virtue of convergence  $\sum_{ n }  a_n  F_- (n) <\infty$ and
condition~\aup, we obtain
 \[
\Sigma_{2,2+} \lez  \sum_{n>  xr/\mu} a_n F_- (n)
 =o(1) =o(\tilde{a}_{x/\mu}).
\]
This proves assertion~(ii).

\smallskip

(iii)~Representation~\eqref{h+} is established using  the same argument as the one
employed to prove part~(ii).  That condition~\eqref{F+o} is necessary and sufficient
for~\eqref{lim_til} follows from Remark~\ref{RemB}.  \hfill$\square$

\begin{Remark}
{\em The case where the tails $F_\pm$ decay as semi-exponential functions can be dealt with in a similar fashion, as  integro-local theorems for  $\pr (S_n \in \D [x))$ that are valid on the whole half-line are available in that case as well (see~\cite{BoMo}, \cite{Mo1}). The  semi-exponential case is ``intermediate'' between the case where the tails are l.r.v.f.'s and the case of exponentially fast decaying tails that we will consider in the next section.
}
\end{Remark}

\section{The case where the moment Cram\'er condition is met}
\label{cramer}

Denote by $\f (\lambda) := \exn e^{\lambda \xi}$ the moment generating function of $\xi$ and set
\[
  \lambda_- := \inf\{\lambda: \f (\lambda) <\infty\}\le 0,\qquad
 \lambda_+ := \sup\{\lambda: \f (\lambda) <\infty\}\ge 0.
  \]
We will assume in this section that the moment Cram\'er condition is satisfied:
\begin{equation*}
 \lambda_+ - \lambda_->0.
\end{equation*}

Denote by  $\lambda_{\rm min}$ the minimum point of $\f (\lambda)$. Clearly, it will also be the minimum point of the (convex) function $L (\lambda) :=\ln\f (\lambda)$ and, since $\mu >0,$ one always has $\lambda_{\rm min}<\lambda_+$. Indeed, if $\lambda_- <0$ then $\lambda_{\rm min}< 0\le \lambda_+,$ while if $\lambda_+ >0$ then $\lambda_{\rm min}\le  0< \lambda_+.$

Since the function $L(\lambda)$ is strictly increasing on  $(\lambda_{\rm min},\lambda_+)$, in that interval  there always exists a unique solution $L^{(-1)}(t)$ to the equation $L(\lambda) =t$ for   $t\in  (L(\lambda_{\rm min}),L(\lambda_+))$. Set
\[
\lambda_q:=  L^{(-1)}(-q)\in (\lambda_{\rm min},\lambda_+), \qquad -q\in (L(\lambda_{\rm min}),L(\lambda_+)).
\]

Now introduce  the ``Cram\'er's transform'' of $\xi$ as a random variable  $\xi^{(\lambda)}$ with the ``tilted'' distribution $\pr (\xi^{(\lambda)} \in dt) =\dfrac{e^{\lambda t}}{\varphi(\lambda)}\, \pr (\xi \in dt).$  It is clear that
\[
L' (\lambda) = \frac{\f' (\lambda)}{\f (\lambda)}= \exn \xi^{(\lambda)}, \qquad \lambda \in (\lambda_-, \lambda_+),
\]
so that $\mu_q := \exn \xi^{(\lambda_q)} \equiv L' (\lambda_q)>0.$

Similarly to notation~\eqref{tia}, for a numerical sequence
$\{b_n\}$ we will denote by   $\tilde b_n$ its ``moving averages'' over intervals of lengths~$d(n)$:
\[
\tilde b_n := \frac1{d (n)}\sum_{n\le k< n+d (n)} b_k , \qquad n\ge 1.
\]

\begin{Theorem}
 \label{thm3}
Assume that  $\lambda_+ - \lambda_->0$ and  $a_n = b_n e^{qn},$ where $\{b_n\}$ satisfies  condition \acu\ with
$\psi (t) =\sqrt{t},$ and  $-q\in (L(\lambda_{\rm min}),L(\lambda_+)).$

{\rm (i)} In the non-lattice case, for any fixed  $0<\D_1<\D_2<\infty,$ the relation
\begin{equation}
 \label{lim_cr}
h (x,\D) \sim \frac{(1-e^{-\lambda_q\D} ) e^{-  \lambda_q x}}{ \mu_q\lambda_q} \, \tilde b_{x/\mu_q}  ,\qquad
x\to\infty ,
\end{equation}
holds uniformly in~$\D\in [\D_1,\D_2].$ Here, if $q=0$ then $\lambda_q=0,$ $\mu_q=\mu$ and, by continuity, the coefficient of~$\tilde b_{x/\mu_q}$ in \eqref{lim_cr}  turns into $\D/\mu,$ as in~\eqref{lim_til}.

{\rm (ii)} In the arithmetic  max-span~1 case, for integer-valued $x\to\infty$ one has
\begin{equation}
\label{lim_ar}
\sum_{n=0}^\infty a_n \pr (S_n=x) \sim \frac{  e^{-  \lambda_q x}}{ \mu_q } \, \tilde b_{x/\mu_q}.
\end{equation}
\end{Theorem}

\proof. We will only present the proof in the non-lattice case; in the arithmetic one, the
argument is even simpler.  Using the standard observation that, for the sum
$S_n^{(\lambda)}:= \xi^{(\lambda)}_1 + \cdots + \xi^{(\lambda)}_n$ of $n$ i.i.d.\
copies of $\xi^{(\lambda)}$, one has $ \pr (S_n^{(\lambda)} \in dt) =\dfrac{e^{\lambda t}}{\varphi^n(\lambda)}\, \pr (S_n  \in dt) $ (this was also used in  Example~1 in~\cite{OmTe}), we obtain
\[
e^{qn}\pr (S_n  \in dt) = e^{(q+L(\lambda))n}e^{-\lambda t} \pr (S_n^{(\lambda)} \in dt).
\]
Therefore,  splitting $\D[x)$ into $K$ subintervals $\D_x[x+ j\D_x),$  $j=0,1,\ldots, K-1,$ of length  $\D_x:=\D /K,$
where $K=K(x)\to\infty$ slowly enough, we have
\begin{align}
h (x, \D) & =
 \sum_n b_n  e^{qn}\pr (S_n \in \D[x))
=
 \sum_n b_n \int_{\D[x)} e^{qn}\pr (S_n \in dt)
 \notag\\
 &  = \sum_n b_n \int_{\D[x)}e^{-\lambda_q t}\pr (S_n^{(\lambda_q)} \in dt)
 \notag\\
 & = \sum_n b_n \sum_{j=0}^{K-1}
  \int_{\D_x[x+ j\D_x)} e^{-\lambda_q t}\pr (S_n^{(\lambda_q)} \in dt)
  \notag\\
 & \sim \sum_n b_n \sum_{j=0}^{K-1}  e^{-\lambda_q (x+ j\D_x) }
    \pr (S_n^{(\lambda_q)} \in \D_x[x+ j\D_x))
    \notag\\
    & =  \sum_{j=0}^{K-1}  e^{-\lambda_q (x+ j\D_x) }
    \sum_n b_n
    \pr (S_n^{(\lambda_q)} \in \D_x[x+ j\D_x)),
    \label{long_f}
\end{align}
where the change of summation order in the last line can be justified using   relation~\eqref{lim_ba} below and condition~\acu.

Now note  that, for  the inner sum   on the right-hand side of~\eqref{long_f}, one has
\begin{equation}
 \label{lim_ba}
 \sum_n b_n
    \pr (S_n^{(\lambda_q)} \in \D_x[x+ j\D_x)) \sim  \frac{\Delta_x}{\mu_q }\,  \tilde b_{( x+ j\D_x)/\mu_q}
  \sim  \frac{\Delta_x}{\mu_q }\,  \tilde b_{x/\mu_q}   , \quad
 x\to\infty,
\end{equation}
uniformly in $j\le K$ (the second equivalence holds since $\{\tilde b_n\}$ is $\psi$-l.c.\ for $\psi(t)=\sqrt{t}$). Indeed, if $\{b_n\}$ satisfies condition \ado, the relation follows from an argument similar to the one proving   assertion~(i) of Theorem~\ref{thm1}, but with  $a_n$ replaced by $b_n$, and $S_n$ replaced by $S_n^{(\lambda_q)}$: the results  of Lemmata~\ref{lem1}, \ref{lem2} and \ref{lem5}(i)
are still applicable in this case, so one only needs to show that,  for the ``left-most sum'' $\Sigma_{2-}$ (now with $m_- = x/(r\mu_q)$), one has  $\Sigma_{2-}=o(\tilde
b_{x/\mu_q})$. This can be easily done using the exponential Chebyshev's inequality,
the obvious representation $\exn e^{\delta\xi^{(\lambda_q)} }=\varphi ( \lambda_q
+\delta)/\varphi ( \lambda_q ) $ (for $\delta\in (0, \lambda_+ - \lambda_q)$) and the
fact that, by Corollary~1 in~\cite{BoBoPap}, one has $\tilde b_n = e^{o(\sqrt{n})}$. Namely,
choosing $\eta:= (r-1)/3>0$ and observing that $\varphi ( \lambda_q +\delta)/\varphi (
\lambda_q )\le e^{\delta \mu_q (1+\eta)} $ for small enough $\delta >0$, we obtain that
\begin{align*}
 \Sigma_{2-}& = \sum_{n<m_-} b_n     \pr (S_n^{(\lambda_q)} \in \D_x[x+ j\D_x))
  \lez
 \sum_{n<m_-} \tilde b_n     \pr (S_n^{(\lambda_q)} \ge x )
 \\
 &\le  \sum_{n<m_-}e^{o (\sqrt{n}) -\delta x}
 \biggl(\frac{\varphi ( \lambda_q +\delta)}{\varphi ( \lambda_q )}\biggr)^n
 \\
 &\le  \sum_{n<m_-} \exp\{o (\sqrt{n}) -\delta x + n  \delta \mu_q (1+\eta)\}
 \\
 &\lez \sum_{n<m_-} \exp\{  -\delta x + n  \delta \mu_q (1+2\eta)\}
 \\
 &\lez   \exp\{  -\delta x + m_-  \delta \mu_q (1+2\eta)\} =e^{- \delta\eta x/r },
\end{align*}
which completes the proof of~\eqref{lim_ba} since, as we said above, $\tilde
b_{x/\mu_q} = e^{o(\sqrt{x})}$. The case where $\{b_n\}$ satisfies condition \aup\ is
dealt with in a similar way.

Thus we showed that the expression in the last line in~\eqref{long_f} is equal to
\[
 (1+o(1) )\frac{\tilde b_{x/\mu_q}}{\mu_q }   \sum_{j=0}^{K-1}  e^{-\lambda_q (x+ j\D_x) } \Delta_x
  \sim
  \frac{\tilde b_{x/\mu_q}}{\mu_q }  \int_{\D [x)} e^{-\lambda_q t} dt
   =
   \frac{(1-e^{-\lambda_q\D} ) e^{-  \lambda_q x}}{ \mu_q\lambda_q} \, \tilde
   b_{x/\mu_q},
\]
which establishes~\eqref{lim_cr}. In the arithmetic case, the proof  proceeds in the
same way. Theorem~\ref{thm3} is proved. \hfill$\square$

\end{document}